\documentclass[12pt,a4paper]{amsart}

\usepackage{amssymb,amsmath,graphics,verbatim}
\usepackage{latexsym}
\usepackage{eucal}
\usepackage{a4wide}

\newtheorem{proposition}{Proposition}[section]
\newtheorem{theorem}{Theorem}[section]
\newtheorem{lemma}[theorem]{Lemma}

\newtheorem{remark}[theorem]{Remark}

\newcommand{\rr}{\mathbb{R}}

\newcommand{\nn}{\mathbb{N}}
\newcommand{\cc}{\mathbb{C}}
\newcommand{\hh}{\mathbb{H}}

\newcommand{\eps}{\epsilon}

\newcommand{\pl}{\partial}
\newcommand{\x}{\times}

\newcommand{\cjd}{\rangle}
\newcommand{\cjg}{\langle}

\newcommand{\ndemi}{\frac{n}{2}}

\def\qed{\hfill$\square$\medskip}

\begin{document}
\title[Spectral characterization of Poincar\'e-Einstein manifolds]{Spectral characterization
of Poincar\'e-Einstein manifolds with infinity of positive Yamabe type}
\author{Colin Guillarmou}
\address{Laboratoire J.A. Dieudonn\'e\\
U.M.R. 6621 CNRS\\
Universit\'e de Nice Sophia-Antipolis\\
Parc Valrose, 06108 Nice\\France}
\email{cguillar@math.unice.fr}

\author{Jie Qing}
\address{Department of Mathematics\\ 
University of California, Santa Cruz\\
CA 95064}
\email{qing@ucsc.edu}

\begin{abstract}
In this paper, we give a sharp spectral characterization of conformally compact Einstein manifolds with conformal infinity of positive
Yamabe type in dimension $n+1>3$. 
More precisely, we prove that the largest real scattering pole of a conformally compact Einstein manifold $(X,g)$ is less than $\ndemi -1$ if
and only if the conformal infinity of $(X,g)$ is of positive Yamabe type. If this positivity is satisfied, we
also show that the Green function of the fractional conformal Laplacian $P(\alpha)$ on the conformal infinity is non-negative for all $\alpha\in [0, 2]$.
\end{abstract}

\maketitle

\begin{section}{Introduction}

Let $\Gamma$ be a convex co-compact group without torsion of orientation preserving 
isometries of the ($n+1$)-dimensional real hyperbolic space $\hh^{n+1}$, and 
let $\Omega(\Gamma)\subset S^n$ the domain of discontinuity of $\Gamma$. 
Then the hyperbolic manifold $X:=\Gamma \backslash\hh^{n+1}$ is conformally
compact with a conformal infinity $M$ which is locally
conformally flat and given by the compact quotient $M=\Gamma\backslash\Omega(\Gamma)$
when we view the elements of $\Gamma$ as M\"obius transformation acting on the closed unit ball of $\rr^{n+1}$. 
In \cite{SY}, Schoen and Yau proved that the Hausdorff dimension $\delta_\Gamma$ of
the limit set $\Lambda(\Gamma)=S^n\setminus\Omega(\Gamma)$ of the group $\Gamma$ is less than
$\ndemi -1$ if the conformal infinity $\Gamma\backslash\Omega(\Gamma)$ is
of positive Yamabe type (we say that a conformal manifold is of
positive Yamabe type if and only if there is a Riemannian metric in
its conformal class whose scalar curvature is positive). Later it was
proved in \cite{Na} that the converse also holds. Sullivan \cite{S} and
Patterson \cite{Pa} also proved that the Poincar\'{e} exponent of the group $\Gamma$ 
is equal to $\delta_\Gamma$. Moreover, in \cite{Pe}, Perry showed that
the largest real scattering pole of
$\Gamma\backslash\hh^{n+1}$ is given by the Poincar\'{e} exponent
$s=\delta(\Gamma)$ (see also \cite{GN} for a characterization of
$\delta(\Gamma)$ in terms of first resonance). Therefore, in this
context, we know that the largest real scattering pole of 
$\Gamma\backslash\hh^{n+1}$ is less than $\ndemi
-1$ if and only if the conformal infinity $\Gamma\backslash\Omega(\Gamma)$ is
of positive Yamabe type. This result which relates the conformal
geometry of the infinity $\Gamma\backslash\Omega(\Gamma)$ to the spectral
property of the conformally compact hyperbolic manifold
$\Gamma\backslash\hh^{n+1}$ has been very intriguing.\\

Later in \cite{Le}, Lee made a clever use of the positive generalized
eigenfunctions to deduce that there is no $L^2$ eigenvalues in $(0,\frac{n^2}{4})$ 
on $(n+1)$-dimensional conformally compact Einstein manifolds $X$
with conformal infinity of nonnegative Yamabe type.  However, the particular case of hyperbolic convex co-compact quotients
mentionned above shows that the absence of $L^2$ eigenvalues does not imply the positivity of the Yamabe type of 
the conformal infinity (the $L^2$-eigenvalues would be scattering poles in $(\ndemi,n)$). 
A simple explicit example is just obtained by taking the quotient of $\hh^3$ by a Fuchsian group $\Gamma$, giving 
rise to an infinite volume hyperbolic cylinder with section the Riemann surface $\Gamma\backslash\hh^2$.
In the introduction of \cite{Le}, Lee asked what would be a sharp spectral
condition for a conformally compact Einstein manifold to have a
conformal infinity of positive Yamabe type. Considering the hyperbolic cases mentionned above, it is then natural to 
ask wether the fact that the largest real scattering pole is less than $\ndemi -1$ on conformally compact Einstein
manifolds is equivalent to positivity of Yamabe type of the conformal infinity. 
In the spirit of the work of Lee \cite{Le}, we are able to give such a spectral characterization of conformally compact Einstein manifolds with 
conformal infinity of positive Yamabe type.\\

Let us first introduce some notations and state our main theorem
precisely. Suppose that $X$ is an $(n+1)$-dimensional smooth manifold with boundary
$\partial X = M$. A metric $g$ on $X$ is said to be
\emph{conformally compact} if, for a smooth defining function $x$ of the boundary
$M$ in $X$, $x^2g$  extends smoothly as a Riemannian metric
to the closure $\bar X$. A conformally compact metric $g$ is complete, has infinite volume, and induces
naturally a conformal class of metrics $[\hat g] = [x^2g|_{TM}]$
(here $x$ ranges over the smooth boudary defining functions). 
As shown in \cite{M1}, the sectional curvature of a conformally compact metric
converges to $-|dx|^2_{x^2g}$ when approaching the boundary $M$. Hence a
metric $g$ on $X$ is naturally said to be \emph{asymptotically hyperbolic} (AH in short) 
if it is conformally compact and the sectional curvatures converge to $-1$ at the boundary. 
A conformally compact Einstein manifold $(X, g)$ is an AH manifold such that $\text{Ric}(g) = -n g$.\\

If $(X, g)$ is an AH manifold, we know (cf. \cite{GL,GR}) that for any representative $\hat g\in [\hat g]$, there is
a unique geodesic defining function $x$ of $\partial X$ associated to the representative $\hat g$ such that the metric $g$
has the geodesic normal form near the boundary
\begin{equation}\label{normal}
g = x^{-2} (dx^2 + g_x)
\end{equation}
where $g_x$ is a one-parameter smooth family of Riemannian metrics on $M$ with
$\hat{g} = \hat g$. 
In Mazzeo \cite{M1} and Mazzeo-Melrose \cite{MM}, it is shown that the spectrum of the (non-negative) Laplacian $\Delta_g$
acting on functions on an AH manifold $(X, g)$ consists of the union of a finite set $\sigma_p(\Delta_g)\subset(0,\frac{n^2}{4})$ of 
$L^2$-eigenvalues, and a half-line of continuous spectrum $[\frac {n^2}4, +\infty)$.  Recently, Joshi-Sa Baretto \cite{JSB} and
Graham-Zworski \cite{GZ} (building on \cite{MM,GuZ,Pe2}), introduced the scattering operators 
$S(s)$ on AH manifolds. For any $s\in \cc$ such that
$$
{\rm Re}(s) \geq \ndemi, \quad s(n-s)\notin \sigma_p(\Delta_g), \quad
s\notin \ndemi + \frac{\nn}{2},
$$
and $f\in C^\infty(\pl X)$, there is a unique solution $v$ to the equation
\begin{equation}\label{1.1}
(\Delta_g - s(n-s))v = 0 
\end{equation}
on $X$ which can be decomposed as follows
\begin{equation}\label{1.2}
v = F x^{n-s} + G x^s , \textrm{ with } F,G\in C^{\infty}(\bar X) \textrm{ and }F|_{\pl X}=f.
\end{equation}
The scattering operator is the linear operator defined on $C^\infty(\pl X)$ by 
\begin{equation}\label{1.3}
S(s)f = G|_{x=0}. 
\end{equation}
If the metric $g_x$ has an even Taylor expansion at $x=0$ in powers of $x$, it is shown in
 \cite{GZ} (see \cite{Gu} for the analysis of the points in
$(n+1)/2-\nn$) that $S(s)$ has a meromorphic continuation to
the complex plane as a family of pseudo-differential operators of
complex order $2s-n$ on $\pl X$. These results extend the
analysis of \cite{GuZ,PP,Pe2} on hyperbolic manifold $\Gamma\backslash \hh^{n+1}$ to the AH class. 
It is proved in \cite{GZ} that $S(s)$ has first order poles 
at $\ndemi +\nn$, the residues of which are the GJMS conformally covariant Laplacian 
on $(\pl X,[\hat{h}])$ constructed in \cite{GJMS} if 
$g$ is asymptotically Einstein. For our purpose, it is
more convenient to consider the renormalized scattering operator 
\begin{equation}\label{1.4}
P(\alpha) := 2^\alpha \frac {\Gamma(\frac \alpha 2)}{\Gamma(-\frac
\alpha 2)} S(\frac {n+\alpha}2). 
\end{equation}
Those $P(\alpha)$ at regular points are conformally covariant
$\alpha$-powers of the Laplacian, they are self-adjoint when
$\alpha$ is real and unitary when ${\rm Re}(\alpha)=0$, moreover
$P(2)$ is the Yamabe operator of the boundary if the bulk space $X$
is (asymptotically) Einstein. We thus call $P(\alpha)$ the 
\emph{fractional conformal Laplacian} for obvious reason. 
The {\it first real scattering pole}
is defined to be the largest real number $s$ such that $\alpha = 2s
- n$ is a pole of $P(\alpha)$.
\begin{theorem}\label{Theorem 1.1} 
Let $(X,  g)$ be a
conformally compact Einstein manifold of dimension $n+1>3$. The first real scattering
pole is less than $\ndemi -1$ if and only if its conformal
infinity $(M, \ [\hat g])$ is of positive Yamabe type.
\end{theorem}

We can also show that
\begin{theorem}\label{Theorem 1.2} 
Let $(X, g)$ be a
conformally compact Einstein manifold of dimension $n+1>3$ with conformal infinity of
positive Yamabe type. Then, for all $\alpha\in (0, 2]$, $P(\alpha)$
satisfies\\
${\rm(a)}$ the first eigenvalue is positive;\\
${\rm(b)}$ $P(\alpha)1$ is positive for any choice of
representative $\hat g$ of the conformal infinity with positive
scalar curvature;\\
${\rm(c)}$ the first eigenspace is generated by a single positive function;\\
${\rm(d)}$ its Green function is nonnegative.
\end{theorem}

\begin{remark}\label{remark1.3}
1) In both cases, it will be clear from the proof that we actually only need to assume that
$$
{\rm Ric} (g) \geq - ng
$$
and that $g_x$ defined in \eqref{normal} has the asymptotic form near the boundary 
\[g_x=\hat{g}- \frac{2x^2}{n-2}\Big({\rm Ric}(\hat g) - \frac {\hat R}{2(n-1)}\hat g\Big)+ O(x^3) \]
where $\hat{g}$ is a metric on $\pl X$, $\text{Ric}(\hat g)$ is its Ricci curvature tensor 
and $\hat R$ its scalar curvature. Metric with this asymptotic `weakly Einstein' structure are discussed by Mazzeo-Pacard 
\cite{MP}.\\
2) Although we do not discuss this here, the smoothness assumption of $g_x$ up to the boundary 
is not necessary, and a restricted smoothness assumption $C^{k,\alpha}(\bar{X})$ for some $k\geq 3$ 
could rather easily be obtained without much modification.\\
3) It is well known that those four properties in Theorem 1.2 all hold
for the conformal Laplacian $P(2)$. However, when
$\alpha\in (0, 2)$, $P(\alpha)$ is a pseudo-differential operator (non-local) 
and it is interesting to see that these four properties continue to hold then.
\end{remark}

Our proof is essentially based on the maximum principle and the existence 
of a positive supersolution for $\Delta_g-s(n-s)$. To construct this supersolution, we use
a special boundary defining function constructed 
by Lee \cite{Le}, which has the advantage of being a positive generalized (non $L^2$) eigenfunction. 
We shall recall some basic facts about conformally compact Einstein manifolds in the next
section. Then in Section 3 we prove Theorem \ref{Theorem 1.1}. Since the proof is rather simple we will carry out some basic calculations 
for the expansions of $F$ for the convenience of the reader. 
Finally in Section 4 we prove Theorem \ref{Theorem 1.2}. 
The crucial issue will be the nonnegativity of the Green function.

\end{section}

\begin{section}{Positive generalized eigenfunctions}

In this Section, we first lay out basic facts about conformally
compact Einstein manifolds,  then we recall the construction
of positive generalized eigenfunctions, following  \cite{Le,CQY,Q}. 
Let $(X,  g)$ be a conformally compact
Einstein manifold with  conformal infinity $(M,  [\hat g])$.
Is is shown in \cite{GL,GR}, that for any representative $\hat g\in [\hat g]$, there is a unique
geodesic defining function $x$ such that the metric has the geodesic normal form
\begin{equation}\label{2.1}
g = x^{-2}(dx^2 + g_x)
\end{equation}
near the boundary. Using this form and considering a Taylor expansion of $g_x$ at $x=0$, Einstein's equations turn into a system which 
can be solved asymptotically (see \cite{FGR,GR}). One
finds that, when $n$ is odd, the metric has an expansion
\begin{equation}\label{2.2}
g_x = \hat g + g^{(2)}x^2 + \ \text{even powers in $x$} \ +
g^{(n-1)}x^{n-1} + g^{(n)}x^n + O(x^{n+1}), 
\end{equation}
and, when $n$ is even,
\begin{equation}\label{2.3}
g_x = \hat g + g^{(2)}x^2 + \ \text{even powers in $x$} \ + hx^n\log
x + g^{(n)}x^n + O(x^{n+2}). 
\end{equation}
When $n$ is odd, $g^{(2i)}$ for $2i<n$ are formally determined by the local
geometry of $(M,  \hat g)$ and $g^{(n)}$ is trace free and
nonlocal. When $n$ is even, $g^{(2i)}$ for $2i<n$, $h$ and the trace
of $g^{(n)}$ are determined by the local geometry of $(M^n, \ \hat
g)$, $h$ is trace free, and trace free part of $g^{(n)}$ is
formally undetermined. Acually, for the purpose of this paper, we only need to assume that
\begin{equation}\label{2.4}
g^{(2)} = - \frac 2{n-2}\Big(\text{Ric}(\hat g) - \frac {\hat
R}{2(n-1)}\hat g\Big), 
\end{equation}
where $\text{Ric}(\hat g)$ is the Ricci curvature tensor of $\hat g$
and $\hat R$ is the scalar curvature of $\hat g$. The following
positive generalized eigenfunction was first constructed and used by
Lee \cite{Le}. Its importance in the results of \cite{Q,CQY} is also worth mentioning. 
From Lemma 5.2 in \cite{Le}, we have
\begin{lemma}\label{Lemma 2.1} 
Let $(X,  g)$ be a conformally
compact Einstein manifold and assume that $\hat g$ is a
representative in $[\hat g]$ of the conformal infinity $(M^n, \
[\hat g])$ and let $x$ be the associated geodesic boundary defining function.
Then there is a unique positive generalized eigenfunction $u$ solving
$$
(\Delta_g + n+1)u = 0
$$
with expansion at the boundary
\begin{equation}\label{2.5}
u = \frac 1x + \frac {\hat R}{4n(n-1)} x + O(x^2). 
\end{equation}
\end{lemma}
The important observation by Lee \cite{Le} (see also
an interesting interpretation of such observation in \cite{Q,CQY})
is that the gradient of $u$ is controlled by $u$: 
\begin{lemma}\label{Lemma 2.2} 
Suppose that, in addition to the assumptions in Lemma \ref{Lemma 2.1}, 
the scalar curvature satisfies $\hat R\geq 0$. Then one has 
\begin{equation}\label{2.6}
|\nabla_g u|_g^2 < u^2 \textrm{ in } X 
\end{equation}
\end{lemma}
\noindent\textsl{Proof}. The proof is done in Proposition 4.2 of \cite{Le}. We repeat it for the the convenience of the reader. 
First the estimate near the boundary 
\begin{equation}\label{2.7}
u^2 - |\nabla_g u|^2 = \frac {\hat R}{n(n-1)} + o(1) 
\end{equation}
follows from the construction of generalized eigenfunctions in Graham-Zworski \cite{GZ}, then an easy computation using 
$\Delta_gu=-(n+1)u$ gives
\begin{equation}\label{2.8}
 \Delta_g (u^2 - |\nabla_g u|_g^2)= 2\cjg({\rm Ric}_g+n) du,du\cjd_g+2\Big|\frac{\Delta_gu}{n+1}g+\nabla_g^2u\Big|^2_g
\end{equation}
which is non-negative if ${\rm Ric}(g)\geq -ng$. From the strong maximum principle,  $u^2-|\nabla u|^2_g$ attains its minimum on $\pl X$ 
and only on $\pl X$, or else is constant. But by \eqref{2.7}, the minimum on $\pl X$ is non-negative, so $|\nabla_gu|_g^2\leq u^2$. 
If $u^2-|\nabla_gu|^2$ is a positive constant, the proof is clearly finished, so it remains to show that  
$u^2$ can not be identically equal to $|\nabla_g u|^2$. If it were the case, an easy computation would give that, for $s>n/2$ and $\phi:=u^{-s}$, 
\[\Delta_g\phi=-s\phi \frac{\Delta_g u}{u}-s(s+1)\phi\frac{|\nabla_g u|^2}{u^2}=s(n-s)\phi\] 
but since clearly $\phi\in L^2$, it contradicts the result of \cite{Le} showing that there is no $L^2$-eigenvalues in $(0,n^2/4)$. 
\qed\\

We remark that for the above two lemmas to hold, we only need to 
assume that $\text{Ric}(g) \geq -n g$ and the expansion \eqref{2.2} and
\eqref{2.3} hold up to second order with $g^{(2)}$ given by \eqref{2.4}.

\end{section}

\begin{section}{Proof of Theorem \ref{Theorem 1.1}}

We present a proof of Theorem \ref{Theorem 1.1} in this section. First we restate
the result of Lee \cite{Le} in terms of scattering pole as follows:

\begin{theorem}\label{Theorem 3.1} 
Let $(X, g)$ be a conformally compact Einstein manifold of dimension $n+1>3$ with conformal infinity of
nonnegative Yamabe type. Then the first scattering pole is less than or equal to $\ndemi$.
\end{theorem}

Here we used the identification of poles of $P(2s-n)$ and poles of the resolvent
$R(s):=(\Delta_g-s(n-s))^{-1}$ in ${\rm Re}(s)>n/2$ (see for instance \cite[Lemma 4.13]{PP}).
Hence to push the first scattering pole down to $\ndemi-1$, we first
show that the scattering operator is regular at $\ndemi$. For this
purpose, we review some of the spectral analysis on AH manifolds. 
By the result of Mazzeo-Melrose \cite{MM,Gu}, the
resolvent of Laplacian $R(s)$ is bounded on $L^2(X)$ for
$$
s\in \cc, \quad {\rm Re}(s)>n/2, \quad
s(n-s)\notin\sigma_{p}(\Delta_g),
$$
and admits a meromorphic continuation to $\cc$ as an operator
mapping the space $\dot{C}^\infty(\bar{X})$ of smooth functions on
$\bar{X}$ vanishing to infinite order at $\partial X$ to the space
$x^s C^{\infty}(\bar{X})$. Moreover the poles of $R(s)$, called {\it
resonances}, are such that the polar part of the Laurent expansion
of $R(s)$ is a finite rank operator. We first observe
\begin{lemma}\label{Lemma 3.2} 
The resolvent $R(s)$ is analytic at $\ndemi$ if and only if there is no function $v\in x^\ndemi
C^\infty(\bar{X})$ such that $(\Delta_g-n^2/4)v=0$.
\end{lemma}
\noindent\textsl{Proof}. It is rather straightforward to see that Lemma 4.9 of
Patterson-Perry \cite{PP} extends to our case, i.e. only a first order
pole is possible for $R(s)$ at $\ndemi$. Indeed, by spectral
theory $\frac n2$ can only be a pole of order at most $2$. If it is
of order $2$, then $n^2/4$ is an $L^2$ eigenvalue for $\Delta_g$ and
the coefficient of order $(s-\frac n2)^{-2}$ is a finite rank
projector on the $L^2$-eigenspace. The analysis of \cite{MM}
(see the proof of Prop 3.3 in \cite{Gu} for details) shows that the
corresponding $L^2$ normalized eigenvectors $(v_k)_{k=1,\dots,K}$ 
would be in $x^\frac n2C^\infty(\bar{X})$, but
to be in $L^2(X)$, this implies actually that  $v_k\in x^{\frac n2+1}C^\infty(\bar{X})$ and by the indicial equation near $\partial
X$,
$$
(\Delta_g-n^2/4)x^jf(y)=-(j-n/2)^2f(y)+O(x^{j+1}), \quad \forall f\in
C^\infty(\partial X)
$$
which implies $v_k=O(x^\infty)$. But Mazzeo's unique continuation
theorem \cite{M2} shows that then $v_k=0$ for all $k$. Then $\frac n2$ can only be a
pole of order $1$ of $R(s)$, in which case the residue of $R(s)$ is finite
rank with range in $\ker(\Delta_g-\frac {n^2}4)\cap x^\frac n2
C^\infty(\bar{X})$.
Conversely assume that $R(s)$ is analytic at $\frac n2$ and that
there is an $u\in x^\frac n2 C^\infty(\bar{X})$ in $\ker(\Delta_g-
n^2/4)$ with leading asymptotic $u\sim x^\frac n2 f_0(y)$ as $x\to 0$.
Then by Graham-Zworski \cite{GZ}, we can construct for a smooth family in
${\rm Re}(s)=n/2$ of solutions $u_s\in
x^{n-s}C^\infty(\bar{X})+x^sC^\infty(\bar{X})$ such that
$u_{n/2}=u$,
$$
(\Delta_g - s(n-s))u_s  =0,
$$
and
$$
u_s = x^{n-s}(f_0+x^2z_s)+ x^{s}(S(s)f_0  + x^2w_s)
$$
where $S(s)$ is the scattering operator, $z_s,w_s$ are smooth
functions on $\bar{X}$ depending smoothly on $s$ on the line
${\rm Re}(s)=n/2$; notice from \cite{GZ} that $u_s$ can be taken to be of
the form
$$
u_s:=x^{n-s}\Phi(s)-R(s)(\Delta_g-s(n-s))(x^{n-s}\Phi(s))
$$
where $\Phi(s)\in C^\infty(\bar{X})$ is smooth in $s$ on the line
${\rm Re}(s)=n/2$ and such that $x^{n/2}\Phi(n/2)=u$ and
$(\Delta_g-s(n-s))(x^{n-s}\Phi(s))=O(x^\infty)$. Since we assumed $R(s)$
analytic at $s=n/2$, then taking the limit as $s\to n/2$ gives
$u=x^\frac n2(f_0+S(n/2)f_0 +O(x^2))$, which implies that
$S(n/2)f_0=0$, but this is not possible since $S(s)$ is unitary on
the line ${\rm Re}(s)=n/2$ (for instance by Section 3 of \cite{GZ}). Thus the
proof is complete.
\qed\\

To show that the resolvent at $\frac n2$ is analytic, we refine slightly Lee's argument.
\begin{lemma}\label{Lemma 3.3} 
Let $(X,  g)$ be a conformally compact Einstein manifold of dimension $n+1>3$, 
with a conformal infinity of nonnegative Yamabe type. Let $k>0$ and consider
\begin{equation}\label{3.1}
\phi = (ku)^{- \frac n2}\log (ku),
\end{equation}
where $u$ is the positive generalized eigenfunction in Lemma \ref{Lemma 2.1} associated with a choice of $\hat g$ of
nonnegative scalar curvature. Then if $k$ is chosen large enough, we have
\begin{equation}\label{3.2}
\Delta_g \phi > \frac {n^2}4 \phi \textrm{ in }X.
\end{equation}
\end{lemma}
\noindent\textsl{Proof}. This is a simple calculation:
$$
\aligned \Delta_g \phi & =- \frac n2\phi \frac {\Delta_g u}u -\frac
n2(\frac n2 + 1)\phi \frac {|\nabla_g u|_g^2}{u^2} +(ku)^{-\frac n2}
((n+1)\frac {|\nabla_g u|_g^2}{u^2} + \frac {\Delta_g u}u)\\& = \frac
{n^2}4 \phi + \frac {n(n+2)}4 \phi (1 - \frac {|\nabla_g u|_g^2}{u^2}) -
(n+1)(ku)^{ -\frac n2} (1 - \frac {|\nabla_g
u|_g^2}{u^2})\\
& = \frac {n^2}4 \phi + (ku)^{-\frac n2}(1 - \frac {|\nabla_g
u|_g^2}{u^2})(\frac {n(n+2)}4 \log (ku) - (n+1))
\\ & > \frac {n^2}4 \phi,
\endaligned
$$
in $X$ provided
$$
\log (ku) > \frac {4(n+1)}{n(n+2)}.
$$
Here we have used \eqref{2.6} of the previous section.
\qed\\

\begin{theorem}\label{Theorem 3.4} 
Let $(X, g)$  be a conformally compact Einstein manifold of dimension $n+1>3$, with  conformal infinity of
nonnegative Yamabe type. Then the resolvent $R(\lambda)$ is regular
at $\frac n2$ and $S(\frac n2) = -{\rm Id}$.
\end{theorem}
\noindent\textsl{Proof}. By Lemma \ref{Lemma 3.2}, we simply need to prove that there is no 
nontrivial function $v$ solving
$$
(\Delta_g - \frac {n^2}4) v = 0 \quad \text{in $X$}
$$
with
$$
v = Fx^\frac n2
$$
for some smooth $F\in C^\infty(\bar{X})$. A straightforward computation gives
\begin{equation}\label{3.3}
\aligned \Delta_g \frac v{\phi} & = \frac {\Delta_g v}{\phi} - 2\cjg \nabla_g v
, \nabla_g \frac 1\phi\cjd_g -  v \frac {\Delta_g \phi}{\phi^2} - v\frac
{2|\nabla_g \phi|_g^2}{\phi^3}\\
& = -(\frac{\Delta_g \phi}\phi - \frac {n^2}4)\frac v\phi +
2\frac{\nabla_g\phi}\phi\cdot \nabla_g \frac v\phi,
\endaligned
\end{equation}
where $\phi$ is defined in \eqref{3.1} in Lemma \ref{Lemma 3.3}. Now, by considering the asymptotic
behaviour of $\phi$ and $v$ at the boundary, we easily see
that
$$
\frac v\phi \to 0
$$
when approaching the boundary. Hence, if there is a negative interior minimum for $v/\phi$ at $p\in X$,
the term $\nabla_g(v/\phi)$ vanishes at $p$ in \eqref{3.3}, but 
since $-((\Delta_g-n^2/4)\phi)/\phi> 0$ in $X$, we deduce that $\Delta_g(v/\phi)$ is positive
near $p$, and this is not possible by applying the strong maximum principle in a small disc around $p$.
We thus have
$$
v\geq 0 \quad \text{on $X$.}
$$
The same argument with an interior maximum shows that $v\leq 0$ and thus $v=0$. 
To see $S(\frac n2) = -\text{Id}$ in this case, 
the proof of Lemma 4.3 in \cite{Gu} can be applied to our case mutatis
mutandis: it shows that the scattering operator at $n/2$ is given by
$$
S(n/2)=-(\text{Id}-2P_0)
$$
where $P_0$ is a projector with respect to $L^2(M,
\text{dvol}_{\hat{g}})$ on the vector space
$$
V:=\{(x^{-\frac n2}u)|_{\partial \bar{X}}; u\in
\text{Range}(\text{Res}_{\frac n2}R(\lambda))\}.
$$
In particular from Lemma 3.2, we obtain $S(\frac n2) = -\text{Id}$.
\qed\\

So far we have improved Theorem \ref{Theorem 3.1} of Lee and obtained that the
first scattering pole is less than $\frac n2$. To push further we
need to show that the scattering operator $S(s)$ for all $s\in
(\frac n2, \frac n2 +1)$ has no kernel. Indeed, from the work of
Joshi-Sa Barreto \cite{JSB} (see \cite{GZ,PP} for the constant curvature
case), we know that
$\tilde{P}(s):=(1+\Delta_{\hat{g}})^{-s/4}P(2s-n)(1+\Delta_{\hat{g}})^{-s/4}$
is a family of bounded Fredholm operators on
$L^2(\partial\bar{X},\text{dvol}_{\hat{g}})$ and the theory of
Gohberg-Sigal \cite{GS} can be used to deduce that, by the meromorphic
functional equation (e.g. see section 3 in \cite{GZ})
$$
S(s) S(n-s) = \text{Id},
$$
the operator $\tilde{P}(2s-n)$ has a pole at $s_0\in \{{\rm Re}(s)\leq n/2\}$ if and only if $\tilde{P}(n-2s_0)$ has a
non-zero kernel, or equivalently $P(2s-n)$ has a pole at $s_0$ if and only if $\tilde{P}(n-2s_0)$ has non-zero kernel.
Thus this corresponds to prove that for $s\in(\frac{n}{2},\frac{n}{2}+1)$, there is no solution to the Poisson equation
\begin{equation}\label{deltags}
(\Delta_g - s(n-s))v = 0 \quad \text{in $X$}
\end{equation}
with
$$
v \in x^{n-s} C^\infty(\bar{X}).
$$
This can be compared to the result of Lee did \cite{Le}: he proved that there is no nontrivial solution 
to the same equation with
$$
v =x^s F \textrm{ for some } F\in C^\infty(\bar{X})
$$
and some $s\in (\frac n2, \frac n2 +1)$. We now define the function
\begin{equation}\label{3.4}
\psi := u^{-(n-s)}. 
\end{equation}
By \eqref{2.5}, we have
\begin{equation}\label{3.5}
\psi = x^{n-s} - \frac {(n-s)\hat R}{4n(n-1)}x^{n+2-s} +
O(x^{n+2-s}). 
\end{equation}
It is also an easy calculation similar to \eqref{3.2} (see also \cite{Le} for the case $\psi=u^{-s}$)
to see that for $s\in(\frac{n}{2},\frac{n}{2}+1)$
\begin{equation}\label{3.6}
\Delta_g \psi > s(n-s)\psi \quad \textrm{ in }X.  
\end{equation}
In order to show that the kernel of $S(s)$ is $0$ for $s\in (\frac
n2, \frac n2 +1)$, we need to find the second term in the expansion
of $F\in C^\infty(\bar{X})$ at the boundary 
(recall $v=x^{n-s}F$ is a solution of \eqref{deltags}). 
This can be found for instance in \cite{GZ}, but
we will give some details for the convenience of
the reader since it is rather straightforward. Recall that, in the
product decomposition $(0,\epsilon)_x\times M$ near the boundary, we
have for any smooth function $f$ defined on $(M, \hat g)$ and any $z\in\rr$
\begin{equation}\label{eqindic}
\begin{split}
\Delta_g (f x^z) =&-\frac {f x^{n+1}}{\sqrt{\det g_x}}
\partial_x (x^{1-n}\sqrt{\det g_x} \partial_x x^z) - \frac
{x^{z+2}}{\sqrt{\det g_x}}\partial_\alpha( \sqrt{\det g_x}
g^{\alpha\beta}_x\partial_\beta f)\\
=& z(n-z) f x^z - \frac z2 f x^{z+1}
\text{Tr}_{\hat g}(\partial_x g_x)+
 x^{z+2} \Delta_{\hat g} f +o(x^{z+2}),
\end{split}
\end{equation}
where $\Delta_{\hat{g}}$ is the Laplacian of $(M, \hat g)$. Hence, since
$$
\text{Tr}_{\hat g}(\partial_x g_x) = - \frac {\hat R}{(n-1)}x + O(x^3)
$$
from \eqref{2.4}, we have
$$
(\Delta_g - s(n-s)) (f x^{n-s}) = (\frac {(n-s)\hat R}{2(n-1)}f
+\Delta_{\hat{g}} f)x^{n-s+2} + o(x^{n-s+2})
$$
and
$$
\aligned (\Delta_g -s(n-s)) (h x^{n-s+2}) & = ((n-s+2)(s-2) -
s(n-s))h x^{n-s+2} + o(s^{n-s+2}) \\ & = -2(n+2-2s) h x^{n-s+2} +
o(x^{n-s+2}). \endaligned
$$
Therefore we have
\begin{equation}\label{3.7}
F = f + \frac 1{2(n+2-2s)}\Big(\frac {(n-s)\hat R f}{2(n-1)} +\Delta_{\hat{g}}
f\Big) x^2 + o(x^2). 
\end{equation}

\begin{lemma}\label{Lemma 3.5} 
Let $(X, g)$ be a conformally compact Einstein manifold of dimension $n+1>3$, with conformal infinity of positive
Yamabe type, and suppose that $h$ is a solution to
$$
S(s) h = 0
$$
on $M$ for some $s\in (\frac n2, \frac n2+1)$. Then $h$ must vanish on
$M$.
\end{lemma}

\noindent\textsl{Proof}. First of all, the statement here is
independent of the choice of representative in $[\hat g]$. We then choose a representative $\hat g$ whose
scalar curvature is positive at every point on $M$. Assume that $h$ is non identically $0$,
we  may assume with no loss of generality that the maximum of $h$ is $1$ and is
achieved at $p_0\in M$. Then we consider
the solution $v$ to the Poisson equation
$$
(\Delta_g - s(n-s)) v = 0
$$
on $X$ with the expansion
$$
v = F x^{n-s} + G x^s
$$
where $F|_{x=0} = h$. Hence, combining \eqref{3.6} and the identity 
\begin{equation}\label{computation}
\Delta_g \frac v{\psi} = -\Big(\frac{\Delta_g \psi}\psi - s(n-s)\Big)\frac v\psi +
2\frac{\nabla_g\psi}\psi\cdot \nabla_g \frac v\psi,
\end{equation}
similar to \eqref{3.3}, we deduce from the maximum principle (exactly like in the proof of Theorem \ref{Theorem 3.4}) that $v/\psi$ can not have an interior 
positive maximum in $X$. The function $v/\psi$ extends continuously to $\bar{X}$ and since its maximum over the boundary is equal to $1$, 
it is clear that $v\leq \psi$ on $X$.
From \eqref{3.7}, we have
\begin{equation}\label{3.8}
v(x, p_0) = x^{n-s} + \frac 1{2(n+2-2s)}\Big( \frac {(n-s)\hat R
}{2(n-1)} + \Delta_{\hat{g}} h (p_0)\Big)x^{n-s+2} +o(x^{n-s+2}). 
\end{equation}
Recall that $p_0$ is a maximum point for $h$ on $M$, which implies that
$\Delta_{\hat{g}} h(p_0) \geq 0$. Comparing  \eqref{3.5} and \eqref{3.8} near $p_0$, 
we obtain a contradiction with the fact that $v\leq \psi$.
\qed\\

It is obvious that Theorem \ref{Theorem 3.4} and Lemma \ref{Lemma 3.5} imply that, for a
conformally compact Einstein manifold with conformal infinity of
positive Yamabe type, the first scattering pole  is less than $\frac
n2 -1$. On the other hand, if we know that the first scattering pole
on an AH manifold is less than $\frac n2 -1$, then we have
$P(0)=\text{Id}$ and so the operator $P(\alpha)$ remains positive for
all $\alpha\in[0,2]$. In particular, the Yamabe operator $P(2)$ is positive and then it is
well known that the conformal infinity is of positive Yamabe type.
This achieves the proof of Theorem \ref{Theorem 1.1}.
\end{section}

\begin{section}{Proof of Theorem \ref{Theorem 1.2}}
Statement (a) in Theorem \ref{Theorem 1.2} is a simple consequence of Theorem \ref{Theorem 1.1}.
Since
$$
P(0) = {\rm Id}
$$
and
$$
P(2) = \Delta_{\hat{g}} + \frac{n-2}{4(n-1)}\hat R
$$
both with positive first eigenvalue, and $P(\alpha)$ for $\alpha\in
(0, 2)$ has no kernel, the first eigenvalue of $P(\alpha)$ has to
be positive for all $\alpha\in (0, 2)$.\\  

Statement (b) follows easily from the arguments used in the proof of Theorem \ref{Theorem 1.1}. 
Let  us give a short proof in the 
\begin{proposition}\label{Proposition 4.4} 
Let $(X, g)$ be a
conformally compact Einstein manifold of dimension $n+1>3$. Suppose that a
representative $\hat g$ of the conformal infinity has positive
scalar curvature on $M$. Then $P_{\hat{g}}(\alpha)1$ is positive for
all $\alpha\in [0, 2]$, where $P_{\hat{g}}$ denotes the operator $P(\alpha)$ defined using $\hat{g}$ for conformal representative in the conformal infinity. 
\end{proposition}
\noindent\textsl{Proof}. Let $v$ be the solution to the Poisson equation
$$
(\Delta_g - s(n-s))v = 0 \quad \text{in $X$}
$$
with 
$$
v = Fx^{n-s} + Gx^s, \quad F,G\in C^\infty(\bar{X})
$$
and expansions 
\begin{equation}\label{4.9}
F = 1 + \frac {(n-s)\hat R}{4(n+2-2s)(n-1)}x^2 + o(x^2),  \quad  G =
S(s)1 + O(x^2),
\end{equation}
where
$$
\alpha = 2s -n \in (0,2).
$$
Let $\psi$ be the positive supersolution of $\Delta-s(n-s)$ defined in \eqref{3.4}, then using \eqref{computation},  
we derive from the maximum principle (exactly like in the proof of Theorem \ref{Theorem 1.1}) that
$$
v < \psi
$$
in $X$. Then, from the expansion \eqref{3.5} and \eqref{4.9}, we first conclude that
$S(s)1 $ has to be non-positive on $M$ for $s\in (\frac
n2, \frac n2+1)$ since  $v-\psi=x^sS(s)1+o(x^s)$. 
Now if $S(s)1$ vanishes at a point $p\in M$, we can consider again the asymptotics \eqref{3.5} and \eqref{4.9} 
along the line $\{y=p;x<\eps\}$ and by positivity of $\hat{R}(p)$ we obtain a contradiction with $v<\psi$ for $x$ small enough.
We thus conclude that $P_{\hat{g}}(\alpha)1 >0$
everywhere on $M$ for all $\alpha \in (0, 2)$. On the other hand,
$P_{\hat{g}}(\alpha)1>0$ holds at $0$ and $2$ obviously. This ends the
proof.
\qed\\

Though, for the differential
operator $P(2)$, the positivity of the first eigenvalue implies the
other three properties due to the maximum principle, it is not so
straightforward for pseudo-differential operators like $P(\alpha)$ for
$\alpha\in(0, 2)$. Of course, the crucial issue is the nonnegativity
of the Green function of the pseudo-differential operators $P(\alpha)$,
or equivalently the non-positivity of the Green function of the
scattering operator $S(s)$ for $s\in (\frac n2, \frac n2+1)$.

By \cite{MM}, outside the diagonal the Schwartz kernel $R(s;m,m')$ of the
resolvent $R(s) = (\Delta_g - s(n-s))^{-1}$ has the regularity
$$
R(s;m,m') \in (xx')^s C^\infty( \bar{X}\times \bar{X}\setminus
\text{diag}_{\bar{X}}).
$$
Consider the Eisenstein function $E(s)\in C^\infty(X\times
\partial\bar{X})$ defined for $s \not= n/2$ and $s$ not a pole of
$R(s)$ by
$$
E(s;m,y'):=(2s-n)[{x'}^{-s}R(s; m, x', y')]_{x'=0}, \quad m\in X,
y'\in\partial\bar{X}
$$
it solves the equation (for all $y'$ fixed in $\partial\bar{X}$)
$$
(\Delta_g-s(n-s))E(s;\cdot,y')=0 \quad \text{ in } X.
$$
From the structure of the resolvent above, we see that for $y'$
fixed in $\partial\bar{X}$, the function $m \to E(s; m, y')$ is in $x^s
C^\infty(\bar{X}\setminus\{y'\})$. Moreover (see \cite{JSB} or \cite{GZ}), the
leading behavior of $E(s;x,y,y')$ as $x\to 0$ (and for $y\not=y'$) is given by
$$
E(s;x,y,y')=x^{s}(S(s;y,y')+O(x))
$$
where $S(s;y,y')$ is the Schwartz kernel of $S(s)$.

For $s\in(n/2,n/2+1)$ such that $S(s)$ is invertible, the Green
kernel of $S(s)$ is given by $S(n-s;y,y')$ by the functional
equation $S(s)S(n-s)=\text{Id}$ (see again \cite{GZ}). The behavior of
$S(n-s;y,y')$ as $y\to y'$ is analyzed in \cite{JSB} (see the Proof of
Theorem \ref{Theorem 1.1} in \cite{JSB} for the computation of the principal symbol of
$S(s)$).

\begin{lemma}\label{Lemma 4.1} 
The leading asymptotic behavior of $S(s;y,y')$ at the diagonal is given by
$$
S(s;y,y')=\frac{\pi^{-\frac n2}\Gamma(s)}{\Gamma(s-\frac
n2)}(d_{\hat g}(y,y'))^{-2s}+O((d_{\hat g}(y,y'))^{-2s+1})
$$
where $d_{\hat g}(\cdot, \cdot)$ denote the distance for the metric
$\hat g$ on $\partial X$. In particular for $s\in (n/2,n/2+1)$, one
has $\Gamma(n/2-s)<0$ so  $S(n-s;y,y')$ tends to $-\infty$ at the
diagonal $\{y=y'\}$ of $\pl X\x\pl X$.
\end{lemma}

With the above understanding of the Green function $S(n-s; y, y')$
of the scattering operator $S(s)$ for $s\in (\frac n2, \frac n2 +1)$
we know that the corresponding Eisenstein function $E(n-s)$ solves
$$
(\Delta_g - s(n-s))E(n-s) = 0 \quad \text{in $X$}
$$
with the expansion
\begin{equation}\label{4.1}
\aligned E(n-s;x,y,y')  = x^{n-s}\Big(& S(n-s;y, y') \\
&+ \frac{x^2}{2(n+2-2s)}(\frac {(n-s)\hat R}{2(n-1)}S(n-s; y,
y') -\Delta_{\hat g} S(n-s; y, y')) \\ & \quad + o(x^2)\Big),
\endaligned
\end{equation}
near the boundary, $y\not=y'$, and where  $y'\in\partial\bar{X}$ is fixed, when
$g$ is at least asymptotically Einstein up to the second order. Let us first deduce the following Lemmas, which will be useful later.
\begin{lemma}\label{Lemma 4.2} 
Let $(X,  g)$ be a conformally
compact Einstein manifold of dimension $n+1>3$ with conformal infinity of positive
Yamabe type. Then the integral kernel $S(n-s; y, y')$ is non-positive for all $y,
y'\in\partial\bar{X}$ and $s\in(\frac n2, \frac n2+1)$.
\end{lemma}
\noindent \textsl{Proof}. The proof runs similarly to the proof of Lemma \ref{Lemma 3.5} except that
$S(n-s; y, y')$ for a fixed $y'\in \partial\bar{X}$ and $s\in (\frac n2,
\frac n2+1)$ is not bounded from  below according to Lemma \ref{Lemma 4.1}.
\qed\\

\begin{lemma}\label{sn-s}
Let $s\in (\ndemi,\ndemi-1)$, then for all fixed $y\in \pl X$, the set $\{y'\in \pl X; S(n-s;y,y')=0\}$ has empty interior in $\pl X$. 
\end{lemma}
\noindent\textrm{Proof}. Assume $S(n-s;y,y')=0$ for some fixed $y\in \pl X$
and $y'$ in an open set $U\subset \pl X$, then by the indicial  equation \eqref{eqindic} 
we deduce easily that $E(n-s;x,y,y')=O(x^\infty)$ for $y\in U$ and by Mazzeo's unique continuation theorem \cite{M2} 
this would imply that $E(n-s;x,y,y')=0$,  which is not possible.
\qed\\

As a consequence of Lemma \ref{Lemma 4.2} we have
\begin{proposition}\label{Proposition 4.3} 
Let $(X, g)$ be a conformally compact Einstein manifold of dimension $n+1>3$, with conformal infinity of
positive Yamabe type. Then, for each $\alpha\in (0, 2)$, the first
eigenspace of $P(\alpha)$ is spanned by a single positive function.
\end{proposition}

\noindent\textsl{Proof}. 
We first produce a positive eigenfunction for $P(\alpha)$ and $\alpha\in (0, 2)$. Since each $P(\alpha)$ for
$\alpha\in (0, 2)$ is invertible and with nonnegative Green
function given by $P(-\alpha)$ (thanks to the functional equation $P(\alpha)P(-\alpha)={\rm Id}$) ,
we look for the eigenfunction of $P(-\alpha)$ as to
maximize 
\begin{equation}\label{4.2}
\frac {\int_M fP(-\alpha)f\text{dvol}_{\hat g}}{\int_M
|f|^2\text{dvol}_{\hat g}}. 
\end{equation}
By Lemma \ref{Lemma 4.2}, we know that
\begin{equation}\label{4.3}
|P(-\alpha)f| \leq P(-\alpha)|f|, 
\end{equation}
hence
$$
\frac {\int_M fP(-\alpha)f\text{dvol}_{\hat g}}{\int_M |f|^2\text{
dvol}_{\hat g}}\leq \frac {\int_M |f|P(-\alpha)|f|\text{dvol}_{\hat
g}}{\int_M |f|^2\text{dvol}_{\hat g}}.
$$
Therefore there is a nonnegative function $f\geq 0$ which is the
first eigenfunction
\begin{equation}\label{4.4}
P(\alpha)f = \lambda(\alpha) f. 
\end{equation}
It is then easily seen that $f$ has to be positive, again due to
Lemma \ref{Lemma 4.2}. Namely, if $f(y)=0$ and $P(-\alpha; y, y')$ is the Green
function of $P(\alpha)$, then,
$$
0 = f(y) = \lambda(\alpha) \int_M P(-\alpha; y,
y')f(y')\text{dvol}_{\hat g}(y'),
$$
which implies $f\equiv 0$. Next we show that, if $h$ is another
eigenfunction of $P(\alpha)$ with eigenvalue $\lambda(\alpha)$, then the ratio $\frac{h}{f}$ has to be a constant
on $M$. We shall use the conformal covariance property of the regularized
scattering operator. Let us denote 
$P_{e^{2\omega}\hat{g}}(\alpha)$ the operator $P(\alpha)$ defined using the conformal representative
 $e^{2\omega}\hat{g}\in [\hat{g}]$ instead of $\hat{g}$, 
 or equivalently using the boundary defining function $e^{\omega}x$. 
 Then  we have by the conformal covariance of  $P(\alpha)$
\begin{equation}\label{4.5}
P_{u^\frac 4{n-\alpha} \hat g}(\alpha) =
u^{-\frac{n+\alpha}{n-\alpha}}P_{\hat g}(\alpha)u, 
\end{equation}
for any positive function $u$ on $M$. Hence
\begin{equation}\label{4.6}
 P_{f^\frac 4{n-\alpha} \hat g}(\alpha) \frac{h}{f}  =
f^{-\frac{n+\alpha}{n-\alpha}}P_{\hat g}(\alpha)h =f^{-\frac{n+\alpha}{n-\alpha}}\lambda(\alpha)h
= f^{-\frac{n+\alpha}{n-\alpha}}(P_{\hat{g}}(\alpha)f)\cdot \frac{h}{f}
 = (P_{f^\frac 4{n-\alpha} \hat g}(\alpha) 1)\cdot \frac{h}{f},
\end{equation}
where 
$$
P_{f^\frac 4{n-\alpha} \hat g}(\alpha) 1 = \lambda(\alpha) f^{-\frac
{2\alpha}{n-\alpha}} >0.
$$
Let $P_{f^\frac 4{n-\alpha}\hat g}(-\alpha; y, y')\geq 0$ be the Green
function of $P_{f^\frac 4{n-\alpha}\hat g}(\alpha)$. Then
\begin{equation}\label{4.7}
\frac{h}{f} (y) = \int_M P_{f^\frac 4{n-\alpha}\hat g}(-\alpha; y, y')((P_{f^\frac 4{n-\alpha}\hat g}(\alpha) 1)\cdot
\frac{h}{f})(y')\text{dvol}_{f^\frac 4{n-\alpha}\hat g}(y'). 
\end{equation} 
Using that
\[
\int_M P_{f^\frac 4{n-\alpha}\hat g}(-\alpha; y, y')(P_{f^\frac 4{n-\alpha}\hat g}(\alpha)1)(y')\text{dvol}_{f^\frac 4{n-\alpha}\hat
g}(y') =1, 
\]
we  deduce from \eqref{4.7}
\[0=\int_M P_{f^\frac 4{n-\alpha}\hat g}(-\alpha; y, y')(P_{f^\frac 4{n-\alpha}\hat g}(\alpha) 1)(y')\Big[
\frac{h}{f}(y)-\frac{h}{f}(y')\Big]\text{dvol}_{f^\frac 4{n-\alpha}\hat g}(y').\]
Since the Green kernel $P_{f^\frac 4{n-\alpha}\hat g}(-\alpha;y,y')$ and $(P_{f^\frac 4{n-\alpha}\hat g}(\alpha)1)(y')$ are respectively 
non-negative and positive by (b) and (d) of Theorem \ref{Theorem 1.1}, we deduce that for all $y\in\pl X$, $\frac{h}{f}(y)=\frac{h}{f}(y')$ for all $y'\not=y$ such that
$P_{f^\frac 4{n-\alpha}\hat g}(-\alpha;y,y')(P_{f^\frac 4{n-\alpha}\hat g}(\alpha)1)(y')\not=0$. But from Lemma \ref{sn-s}, we know that 
for each $y$, this set is dense in $\pl X$. By continuity of $h$ and $f$ (which follows from ellipticity of $P(\alpha)$), 
we can conclude that $h=f$.
Thus the proof is complete.
\qed\\

\end{section}

\textbf{Acknowledgement}. Both authors thank the Institute for Advanced Study in Princeton where this work was done, in particular C.G.  was partially supported there by NSF fellowship No. DMS-0635607.  We are also grateful to the anonymous referee for 
his careful reading.

\end{document}